\documentstyle{amsppt}
\magnification=1200

\def\R{{\bold R_{\max}}}

\def\a-{{\it a-}}
\def\o-{{\it o-}}
\def\wo{$wo$}
\def\G-{{\it G-}}
\def\Inf{$\wedge$}

\def\0{{\bold 0}}
\def\1{{\bold 1}}

\def\< {\prec }
\def\LE {\preceq }
\def\> {\succ }
\def\GE {\succeq }

\def \inf{{\text {inf}}}
\def \sup{{\text {sup}}}

\NoRunningHeads

\document

{\eightpoint 
\noindent Doklady Mathematics, Vol. 62, No. 2, 2000, pp. 169--171.
Translated from Doklady Akademii Nauk, Vol. 374, No. 1, 2000,
pp. 26--28.

\noindent Original Russian Text Copyright \@ 2000 by Shpiz.

\noindent English Translation Copyright \@ 2000 by MAIK "Nauka/Interperiodika"
(Russia)}

\bigskip

\centerline{MATHEMATICS}

\bigskip

\centerline{\bf A Theorem on Eigenvectors in Idempotent Spaces}

\medskip

\centerline{G. B. Shpiz}

\centerline {Presented by Academician V.P. Maslov. March 25, 2000}

\centerline{Received March 31, 2000}

\bigskip\bigskip

This paper suggests an algebraic version of the theorem on the
existence of eigenvectors for linear operators in abstract idempotent
spaces (Theorem 3). Earlier, the theorem on the existence of
eigenvectors was only known for the cases of a free finite-dimensional
semimodule [3] and for compact operators in semimodules of real
continuous functions [4]. 

{\bf 1}. We use the terminology from [6--8]. Recall that an {\it idempotent 
semigroup} (IS) is an additive semigroup with the commutative addition 
$\oplus $ such that $x\oplus x=x$ for any element $x$. An arbitrary IS 
can be treated as an ordered set with the partial order defined as
follows: $x\LE y$ if an only if $x\oplus y=y$. It is easy to see that 
this order is well-defined and $x\oplus y=\sup (x,y)$. For an arbitrary
subset $X$ of an idempotent semigroup, we put $\oplus X=\sup (X)$ 
and $\wedge X=\inf (X)$ if the corresponding right-hand sides exist.
An IS is called $b$-{\it complete} (or {\it boundedly complete}) if
 any of its subsets
bounded from above (including the empty subset) has the least upper
bound. In particular, any $b$-complete IS contains zero (denoted by $\0 $), 
which coincides with $\oplus \varnothing$, where $\varnothing$
is the empty set. A homomorphism of $b$-complete idempotent semigroups
is called a $b$-{\it homomorphism} if $g(\oplus X)=\oplus g(X)$ for any subset
$X$ bounded from above. An {\it idempotent semiring} (ISR) is an IS endowed with
an associative multiplication $\odot$ with an identity element (denoted by
$\1$) such that both distributivity laws are valid. An ISR is called a
{\it semifield} if any of its nonzero elements has an inverse. An
idempotent semifield is called $b$-{\it complete} if it is $b$-complete as an IS. 
The set of all nonzero elements of a $b$-complete semifield is a complete 
group (in the terminology of [9]). The converse assertion is also valid:
adding $\0$ to any complete ordered group and defining addition as $\sup$, 
we obtain a $b$-complete semifield. The semifield obtained by this method 
from the additive group of all real numbers is denoted by $\R$. The theory of 
ordered groups [9] implies that any $b$-complete semifield is commutative. In 
any $b$-complete semifield, the generalized distributive laws
$$
a\odot (\oplus X)=\oplus (a\odot X),\quad 
a\odot (\wedge X)=\wedge (a\odot X)
$$
are valid; here $a$ is an element of the field and $X$ is a nonempty 
bounded subset. 

An {\it idempotent semimodule} over an idempotent semiring $K$ is an idempotent
semigroup $V$ endowed with a multiplication $\odot$ by elements of $K$
such that, for any $a,b\in K$ and $x,y\in V$, the usual laws
$$
\align
a\odot (b\odot x)&=(a\odot b)\odot x, \\
(a\oplus b)\odot x&=a\odot x\oplus b\odot x,\\
a\odot (x\oplus y)&=a\odot x\oplus a\odot y,  \\ 
\0 \odot x&=\0
\endalign
$$
are valid. An idempotent semimodule over an idempotent semifield is called
an {\it idempotent space}. An idempotent $b$-complete space $V$ 
over a $b$-complete semifield
$K$ is called an {\it idempotent} $b$-{\it space} if, for
 any nonempty bounded subset 
$Q\subset K$ and any $x\in V$, the relations
$$
(\oplus Q)\odot x=\oplus (Q\odot x),\quad 
(\wedge Q)\odot x=\wedge (Q\odot x)
$$
hold. A homomorphism $g: V\to W$ of $b$-spaces is called a $b$-{\it homomorphism},
or a $b$-{\it linear} mapping, if $g(\oplus X)=\oplus g(X)$ for any bounded
subset $X\subset V$. More general definitions 
(for spaces which may not be $b$-complete) see in [7].
Homomorphisms taking values in $K$ (treated as a
semimodule over itself) are called {\it linear functionals}. A subset of an
idempotent space is called a {\it subspace} if it is closed with respect to
addition and multiplication by coefficients. A subspace in a $b$-space is 
called a $b$-{\it closed subspace} if it is closed with respect to 
summation over arbitrary bounded (in V) subsets. This subspace has a natural
structure of $b$-space; it is also a $b$-subspace in V in the sense of [7].

For an arbitrary set $X$ and an idempotent space $V$ over a semifield $K$,
we use $B(X,V)$ to denote the semimodule of all bounded mappings from $X$ 
into $V$ with pointwise operations. If $V$ is an idempotent $b$-space,
then $B(X,V)$ is a $b$-space. A mapping $f$ from a topological space $X$ 
into an ordered set $V$ is called {\it upper semicontinuous} if, for any
$b\in V$, the set  $\{x\in X | f(x)\GE b\}$ is closed in $X$, see [7]. 
In the case
where $V$ is the set of real numbers, this definition coincides with the usual
definition of upper semicontinuity of a real function. The set of all
bounded upper semicontinuous mappings from $X$ to $V$ is denoted by 
$USC(X,V)$. If $V$ is a boundedly complete lattice, then $USC(X,V)$ is
also a boundedly complete lattice with respect to the pointwise order.
If V is an idempotent $b$-space, then $USC(X,V)$ is also a $b$-space
with respect to the operations $f\oplus g=\sup (f,g)$ and 
$(k\odot f)(x)=k\odot f(x)$. 

{\bf 2.} In what follows, unless otherwise specified, the symbol $K$
stands for a $b$-complete idempotent semifield and all idempotent spaces
are over  $K$. 

A subset $M$ of idempotent $b$-space $V$ is called \wo-{\it closed} if 
$\wedge X\in M$ and $\oplus X\in M $ for any linearly ordered subset 
$X\nomathbreak\subset\nomathbreak M$ in $V$. A nondecreasing mapping
$f: V\to W$ of $b$-spaces is called \wo-continuous if  
$f(\oplus X)=\oplus f(X)$ and $f(\wedge X)= \wedge f(X)$ for any bounded
linearly ordered subset $X\subset V$. Note that an arbitrary isomorphism of
ordered sets is \wo-continuous. It can be shown that the notions of 
\wo-closedness and \wo-continuity coincide with the closedness and continuity 
with respect to some $T1$ topology defined in an intrinsic way in terms of the
order.

\proclaim{Proposition 1} Suppose that $V$ is an idempotent b-space and 
$W$ is a \wo-closed subsemigroup of $V$. Then $\oplus \nomathbreak X
\in \nomathbreak W$ for any subset $X\subset W$ bounded in $V$. In
particular, each \wo-closed subspace is a b-closed subspace.\endproclaim

An element $x$ of an idempotent space $V$ is called {\it Archimedean} if, for any
$y\in V$, there exists a coefficient $\lambda \in K$ such that  
$\lambda \odot x\GE y$. For an Archimedean element $x\in V$, the
formula $x^*(y)=\wedge \{k\in K | k\odot x\GE y\}$ defines a mapping 
$x^*: V\to K$. If $V$ is an idempotent $b$-space, then $x^*$ is a
$b$-linear functional and $x^*(y)\odot x\GE y$ for any $y\in V$ [6]. 
We say that an Archimedean element $x\in\nomathbreak V$ is \wo-{\it continuous}
if the functional $x^* $ is \wo-continuous, and that an idempotent 
$b$-space $V$ is {\it Archimedean} if $V$ contains a \wo-continuous Archimedean
element.

\proclaim{Proposition 2} If $X$ is a compact topological space, then 
USC(X,K) is an Archi\-medean space and the function $\bold e$ identically
equal to $\1$ is a \wo-continuous Archime\-dean element.\endproclaim

Note that $\bold {e^*}(f)=\sup \{f(x) | x\in X\}$.

\proclaim{Theorem 1} Any \wo-closed subspace of an Archimedean space
is an Archimedean space. Any linearly ordered (with respect to the
inclusion) family of nonzero \wo-clos\-ed subspaces of an Archimedean
space $V$ has a nonzero intersection.\endproclaim

{\bf 3.} An Archimedean idempotent space $V$ is called {\it irreducible} with
respect to a set $G$ of arbitrary mappings of $V$ into itself if any
nonzero \wo-closed $G$-invariant subspase $W$ of $V$ coincides with $V$.

\proclaim{Theorem 2} Let $V$ be an Archimedean space and $G$ be an
arbitrary set of mappings of $V$ into itself. Then $V$ contains a 
\wo-closed $G$-invariant irreducible Archimedean subspace.\endproclaim

Let $V$ be a semimodule over $K$. An {\it eigenvector} of a mapping $g: V\to V$ 
corresponding to an {\it eigenvalue} $\lambda \in K$ is, by definition, a
nonzero vector $x\in V$ such that $g(v)=\lambda \odot x$. A semiring
$K$ is called {\it algebraically closed} [3] if, for any elements $x\in K$ 
and any positive $n$, there exists an element $y\in K$ such that $y^n=x$. 
For example, $\R$ is a $b$-complete algebraically closed semifield.

\proclaim{Theorem 3} An arbitrary $b$-linear mapping of an Archimedean
idempotent space over an algebraically closed semifield into itself
has an eigenvector.\endproclaim

\proclaim{Proposition 3} Suppose that $V$ is an Archimedean semimodule 
over $\R$, $g$ is a homomorphism (or even an arbitrary homogeneous mapping) 
of the semimodule $V$ into itself, and $x,y\in V$ are eigenvectors of
$g$ corresponding to eigenvalues $p,q\in \R$, respectively. Then  $(i)$ $x\LE y$
implies $p\LE q$ and $(ii)$ all eigenvalues corresponding to Archimedean
eigenvectors coincide.\endproclaim

{\bf 4.} Let $V$ be a $b$-space. A subset $W\subset V$ is called a 
\Inf-{\it subspace} if it is closed with respect to multiplication by scalars
and taking greatest lower bounds of nonempty subsets. By this definition, any
such $W$ is a boundedly complete lattice with respect to other inherited
from the ambient space. Therefore, any \Inf-subspace $W\subset V$ can
be treated as a semimodule with respect to the inherited multiplication
by scalars and operations $x\oplus _W y=\sup (x,y)$, where $\sup$ 
is over $W$. In what  follows, all \Inf-subspaces are considered as 
semimodules with respect to these operations. The definitions
immediately imply that any \Inf-subspace of a $b$-space is a $b$-space. 
It is easy to show that $USC(X,V)$ is a \Inf-subspace in $B(X,V)$ 
for any $b$-space $V$ and any topological space $X$. 

\proclaim{Theorem 4} If $V$ is an Archimedean $b$-space and $x\in V$ 
is a \wo-continuous Archi\-medean element, then any \Inf-subspace $W$ 
of $V$ containing $x$ is an Archimedean $b$-space.\endproclaim

\proclaim{Theorem 5} An idempotent $b$-space $V$ over an algebraically
closed $b$-complete semifield $K$ is Archimedean if and only if there exists 
a space of the form USC{\rm(}X,K{\rm)}, where $X$ is a compact topological space,
such that $V$ is isomorphic to its \Inf-subspace containing constants.
\endproclaim

\bigskip

\centerline{ACKNOWLEDGMENTS}

The author thanks G.L. Litvinov and V.P. Maslov for their help and 
attention.

This work was financially supported by the Dutch Organization
for Scientific Researches (N.W.O.).

\bigskip

\centerline{REFERENCES}

1. Birkhoff, G., {\it Teoriya reshetok} (Lattice Theory), Moscow: Nauka,
1984.
 
2. Dudnikov, P.S. and Samborskii, S.N., {\it Izv. Akad. Nauk SSSR, Ser.
Mat.}, 1991, vol. 55, no. 1, pp. 93--109.
 
3. Dudnikov, P.S. and Samborskii, S.N., {\it Idempotent Analysis}, Adv. Sov.
Math., vol. 13, Providence: Am. Math. Soc., 1992, pp. 65--85.
 
4. Maslov, V.P. and Kolokoltsov, V.N., {\it Idempotentnyi analiz i ego
primeneniya v optimal'nom upravlenii} (Idempotent Analysis and Its
Applications to Optimal Control), Moscow: Nauka, 1994.
 
5. Kolokoltsov, V.N. and Maslov, V.P., {\it Idempotent Analysis and
Applications}, Dordrecht: Kluwer, 1997.
 
6. Litvinov, G.L., Maslov, V.P. and Shpiz, G.B., {\it Dokl. Akad. Nauk},
1998, vol. 363, no. 3, pp. 298--300.
 
7. Litvinov, G.L., Maslov, V.P. and Shpiz, G.B., Idempotent Functional
Analysis: Algebraic Approach, {\it Preprint of Int. Center Sophus Lie},
Moscow, 1998 (http://www.chat.ru/$\sim$sophus\_lie/funana.zip).
 
8. Litvinov, G.L., Maslov, V.P. and Shpiz, G.B., {\it Mat. Zametki},
1999, vol. 65, no. 4, pp. 572--585.
 
9. Fuchs, L., {\it Partially Ordered Algebraic Systems}, Oxford, 1963.
Translated under the title {\it Chastichno-uporyadochennye 
algebraicheskie sistemy}, Moscow: Mir, 1965.

\enddocument